\documentclass[conference]{IEEEtran}

\usepackage{cite}

\ifCLASSINFOpdf
\else
\fi
\usepackage[cmex10]{amsmath}
\usepackage{amsthm}
\usepackage{amssymb}
\usepackage[all,2cell]{xy}\UseAllTwocells\SilentMatrices

\newtheorem{theorem}{Theorem}
\newtheorem{proposition}[theorem]{Proposition}

\newtheorem{lemma}[theorem]{Lemma}
\theoremstyle{definition}
\newtheorem{definition}[theorem]{Definition}
\newtheorem{remark}[theorem]{Remark}
\newtheorem{example}[theorem]{Example}

\begin{document}
\title{Universal Factorizations of Quasiperiodic Functions}

\author{\IEEEauthorblockN{Michael Robinson}
\IEEEauthorblockA{Department of Mathematics and Statistics\\
American University\\
4400 Massachusetts Ave NW\\
Washington, DC 20016\\
Email: michaelr@american.edu}}

\maketitle

\begin{abstract}
Chirped sinosoids and interferometric phase plots are functions that are not periodic, but are the composition of a smooth function and a periodic function.  These functions functions factor into a pair of maps: from their domain to a circle, and from a circle to their codomain.  One can easily imagine replacing the circle with other phase spaces to obtain a general quasiperiodic function.  This paper shows that under appropriate restrictions, each quasiperiodic function has a unique universal factorization.  Quasiperiodic functions can therefore be classified based on their phase space and the phase function mapping into it.
\end{abstract}


\IEEEpeerreviewmaketitle

\section{Introduction}

Chirped sinosoids and interferometric phase plots are functions that are not periodic, but are the composition of a smooth function and a periodic function.  These functions functions factor into a pair of maps: from their domain to a circle, and from a circle to their codomain.  Functions that are the composition of a smooth map and a periodic function can be rather different from almost-periodic functions, since chirps do not lie in the closure of the space of trigonometric polynomials with any reasonable norm.  Additionally, chirps are not periodic to first order, so they are not amenable to asymptotic approaches like the famed WKBJ method \cite{Kevorkian_1996}.

One can easily imagine replacing the circle with other \emph{phase spaces}.  This paper shows that under appropriate restrictions, each such \emph{quasiperiodic} function has a universal factorization.  These quasiperiodic functions can then be classified canonically based on their phase space.  We expect that this canonical representation could serve as the starting point for practical filtering algorithms for quasiperiodic signals.  

\subsection{Contributions}
This paper presents a general definition of a quasiperiodic function and factorizations of it that ensure the phase space is a manifold.  Given this definition, it proves that unique universal factorizations of quasiperiodic functions exist.  To connect the theory of quasiperiodicity to periodic functions, it completely characterizes universal $S^1$-quasiperiodic factorizations of periodic functions.  This paper ends by presenting an initial algorithmic framework for producing factorizations of quasiperiodic functions.

\subsection{Historical context}

The problem of recovering the phase of an analytic signal from its amplitude is called ``phase retrieval'', and has occupied the attention of researchers at least since the 1950s in the context of x-ray crystallography \cite{Calderon_1952}.  Work continues even recently (for instance \cite{Iwen_2015}), though most of the results relevant to this article were obtained in the 1980s \cite{Fienup_1982}.  Numerous methods, such as those based on iteration \cite{Gerchberg_1972}, asymptotics \cite{Kevorkian_1996}, or gradient search methods \cite{Reed_1983}, have been developed to address the problem, though they all rely heavily on the assumption that the phase lies on a circle.  Although phase ambiguities are prevalent in time series data, they are quite rare in the case of a two-dimensional domain \cite{Barakat_1984}.

The study of non-circular phase spaces has been pursued most aggressively by researchers in dynamical systems.  Most of their approaches are based on the famous embedding theorems of Whitney \cite{Whitney_1936}.  In particular, the Takens delay-embedding theorem (see for instance \cite{Takens_1981,Sauer_1991}) can obtain topologically-accurate models of phase spaces that are suitable for factorizations of quasiperiodic functions.

Our approach is inspired by the work by DeSilva \emph{et al.} \cite{DeSilva_2011}, which uses persistent cohomology to construct circular coordinates.  They use a delay embedding approach to obtain a model of the underlying phase space, followed by persistent cohomology to identify nontrivial 1-cycles, which they smooth using harmonic 1-forms.  The class of functions they consider are precisely those that factor into a diffeomorphism and a periodic function.  We generalize beyond their framework in a complementary way, by exploring more general factorizations and justify that their delay-embedding approach still works in our case (though not necessarily the harmonic smoothing).

\section{Motivating examples}

Every smooth periodic function $f:\mathbb{R}\to N$ from the real line to a manifold $N$ can be written as the composition $F \circ \phi$, where $\phi: \mathbb{R} \to S^1$ is a covering map.  The function $F:S^1 \to N$ encodes the behavior of $f$ on a fundamental domain.  The function $\phi$ extends $F$ periodically from $S^1$ to the whole real line.  Since $\phi$ is a covering map, it is a surjection and a local homeomorphism.  Thinking about the topological spaces involved, $S^1$ can be written as the quotient $\mathbb{R}/\phi$, which we call the \emph{phase space} for the factorization $F \circ \phi$.  

There is, however, a different quotient that may be taken as well, namely $\mathbb{R}/f$.  This quotient serves as the phase space for a different factorization $f=G \circ \pi$, where $\pi:\mathbb{R}\to\mathbb{R}/f$ is the canonical projection.  In this factorization, $\pi$ is surjective and $G$ is injective.

Since $f=F\circ\phi$, any pair of points $x,y \in \mathbb{R}$ with $\pi(x)=\pi(y)$ will also satisfy $f(x)=f(y)$.  Therefore, $\mathbb{R}/f$ might appear to be a better object to work with than $\mathbb{R}/\phi$ -- since it is the minimal specification of the function $f$ -- but $\mathbb{R}/f$ isn't usually a manifold.  Wherever $f$ has a critical point, the corresponding point in $\mathbb{R}/f$ will usually not have a Euclidean neighborhood.

\begin{example}
Consider the function $f:\mathbb{R}\to\mathbb{R}$ given by $f(x)=\sin x$.  Then $\mathbb{R}/f$ is the closed interval $[-\pi/2,\pi/2]$, which is not a manifold.  (It is a manifold with boundary, however.)
\end{example}


\section{Quasiperiodicity}
In order to ensure that the phase space of a factorization of a periodic function is a manifold, we need to constrain the projection to the phase space.  The lesson from differential topology is that this projection should be a surjective submersion (see \cite[Prop 7.16-7.19]{Lee_2003} for instance).

\begin{definition}
Suppose $M, N, C$ are smooth manifolds.  A smooth function $u:M \to N$ is called \emph{$C$-quasiperiodic} if it has a factorization $u=U \circ \phi$ where $\phi:M \to C$ is a surjective submersion.  Such a factorization is called a \emph{$C$-quasiperiodic factorization}.  We call $C$ the \emph{phase space} of the factorization.
\end{definition}

Quasiperiodic factorizations of a given function are not unique, even for a given phase space.

\begin{example}
Consider the periodic function $u(x)=\sin 2\pi x$ as a function $\mathbb{R}\to \mathbb{R}$.  One quasiperiodic factorization uses $\mathbb{R}$ as the phase space and $\phi$ as the identity function.  The function $u$ also has many $S^1$-quasiperiodic factorizations.  Using the phase space $C=S^1=\mathbb{R}/\mathbb{Z}$, one factorization has $\phi(x) = x/2 \mod 1$ and $U(\theta)=\sin(4\pi \theta)$.  Another factorization is $\phi'(x) = x/3 \mod 1$ followed by $U(\theta)=\sin( 6\pi \theta)$.
\end{example}

Fortunately, there is a canonical notion of a quasiperiodic factorization, namely one that satisfies a universal property.  This neatly balances the need for a minimal specification with the requirement that the phase space be a manifold.

\begin{definition}
\label{df:universality}
A $C$-quasiperiodic factorization $u=U \circ \phi$ of a smooth function $M\to N$ is \emph{universal} if for any other quasiperiodic factorization $u=U' \circ \phi'$ with $\phi' :M\to C'$, there exists a smooth function $c:C'\to C$ such that the diagram below commutes:
\begin{equation*}
\xymatrix{
M \ar[r]^{\phi'} \ar[d]_\phi & C' \ar[d]^{U'} \ar[dl]^c \\
C \ar[r]_U & N \\
}
\end{equation*}
\end{definition}

A universal quasiperiodic factorization has the minimal phase space among quasiperiodic factorizations.  It therefore generalizes the notion of a fundamental domain for a periodic function.  This connection is made precise in Theorem \ref{thm:s1_universality}.  

\begin{theorem}
\label{thm:universality_exists}
Every quasiperiodic function has a universal quasiperiodic factorization, which is unique up to a diffeomorphism applied to the phase space.
\end{theorem}

The proof of this theorem is given in Section \ref{sec:general_universality}, and a complete characterization of universal $S^1$-quasiperiodic factorizations of periodic functions is given in Section \ref{sec:s1_universality}.  

\begin{remark}
For a function $u:M\to N$, the quotient $M/u$ also satisfies a universal property like the one described in Definition \ref{df:universality}.  However, as noted in the previous section, $M/u$ is generally not a manifold.
\end{remark}

Because of their uniqueness, universal quasiperiodic factorizations govern the structure of all quasiperiodic factorizations.  Therefore if we find a quasiperiodic factorization by any means, we know its structure is constrained from the outset: it must have the universal phase space as a submanifold and must descend to this submanifold as a quotient.  

In some cases, universality of a quasiperiodic factorization is easy to detect as the next Proposition shows.

\begin{proposition}
Suppose $u:M\to N$ is a constant function.  Then $u:M\to * \to N$ is its universal quasiperiodic factorization, where $*$ is the space with one point.  On the other hand, if $v=V \circ \phi$ is a quasiperiodic factorization where $V$ is injective, then it is a universal quasiperiodic factorization.
\end{proposition}
\begin{IEEEproof}
For the first statement, consider any other factorization $u:M\to C\to N$.  Construct the unique function $c:C\to *$ and observe that the diagram shown below left commutes:
\begin{equation*}
\xymatrix{
M \ar[d] \ar[r] & {*} \ar[d] & M  \ar[d]_{\phi'}\ar[r]^\phi & C \ar[d]^V \\
C \ar[r] \ar[ru]^c & N & C' \ar[r]_{V'} \ar[ru]^c & N \\
}
\end{equation*}

For the second statement, consider the diagram above right and suppose that $v=V' \circ \phi'$ is another quasiperiodic factorization.  Since $V$ is injective, it has a partially defined inverse $V^{-1}:v(M) \to C$, which means that $c:C'\to C$ given by $c=V^{-1} \circ V'$ completes the commutative diagram.
\end{IEEEproof}

\section{Universality for periodic functions}
\label{sec:s1_universality}

This section presents a complete characterization of the universal $S^1$-quasiperiodic factorizations of periodic functions, thereby establishing a link between the theory of quasiperiodic functions and the traditional theory of periodic functions.  The phase space of such a universal quasiperiodic factorization is precisely its fundamental domain.

\begin{theorem}
\label{thm:s1_universality}
Suppose $u:\mathbb{R}\to N$ is a smooth, non-constant $T$-periodic function.  Then
\begin{enumerate}
\item $u$ is $S^1$-quasiperiodic and has a phase function $\phi:\mathbb{R} \to S^1 = \mathbb{R}/\mathbb{Z}$ given by $\phi(x)=x/T \mod 1$.
\item The factorization is universal unless there is an $S^1$-quasiperiodic factorization of $u$ whose phase function is non-injective on an interval of the form $[a,a+T)$.  
\end{enumerate}
\end{theorem}

\begin{example}
As an example of a universal quasiperiodic factorization of a periodic function, consider the $S^1$-quasiperiodic function defined on $[0,2]$ given by the factorization $f=U\circ \phi$ where
\begin{equation*}
\phi(x) = x/2 \mod 1
\end{equation*}
and
\begin{equation*}
U(\theta)=\begin{cases}
\sin 4\pi \theta & \text{if } 0 \le \theta < 1/2\\
\frac{1}{2} \sin 8\pi \theta & \text{if } 1/2 \le \theta < 1.\\
\end{cases}
\end{equation*}
This is an example of a universal quasiperiodic factorization.  Notice that any factorization that has a phase function with a period smaller than $2$ would be unable to correctly reproduce the change in amplitude of $u$.
\end{example}

\begin{example}
The smoothness of $u$ is a necessary condition for Theorem \ref{thm:s1_universality}.  Consider the case of the function given by the factorization $g=U \circ \phi$ where
\begin{equation*}
\phi(x) = x/2 \mod 1 
\end{equation*}
and 
\begin{equation*}
U(x) = \begin{cases}
\sin 4\pi \theta & \text{if } 0 \le \theta < 1/2\\
\sin 8\pi \theta & \text{if } 1/2 \le \theta < 1.\\
\end{cases}
\end{equation*}
almost as in the previous example, though $g$ is not smooth.  However, there is also a new factorization $g=U' \circ \phi'$
\begin{equation*}
\phi'(x) = \begin{cases}
x & \text{if } 0 \le x \mod 2 < 1 \\
2x-2 & \text{if } 1 \le x \mod 2 < 2,\\
\end{cases}
\end{equation*}
with 
\begin{equation*}
U'(\theta) = \sin 2\pi \theta.
\end{equation*}
Notice that $\phi'$ has fundamental period 2, but that there is no map $c$ that makes the appropriate diagram commute -- contradicting the universality of $g = U \circ \phi$.
\end{example}

\begin{proposition}
Any function that is the composition of a diffeomorphism $D:\mathbb{R}\to\mathbb{R}$ and a periodic function $p:\mathbb{R} \to N$ is $S^1$-quasiperiodic.
\end{proposition}
\begin{IEEEproof}
The proof is a simple commutative diagram:
\begin{equation*}
\xymatrix{
\mathbb{R} \ar[r]^D \ar@/_1pc/[rr]_{\phi'=\phi\circ D} & \mathbb{R} \ar@/^1.5pc/[rr]^{p} \ar[r]^\phi & S^1 \ar[r] & N.\\
}
\end{equation*}
\end{IEEEproof}

Lest one think that $S^1$-quasiperiodic functions are all essentially just periodic functions with distorted timescales, the converse is not true.

\begin{example}
\label{eg:not_even_periodic}
Consider the function $u:\mathbb{R}\to S^1=\mathbb{R}/\mathbb{Z}$ given by $u(x)=\arctan x \mod 1$.  This function is $S^1$-quasiperiodic -- trivially -- since $\arctan x$ is a submersion and has a range larger than $[0,1)$.  This is completely different than a periodic function since a preimage of $u$ contains at most $3$ points -- not an infinite set.
\end{example}

\begin{IEEEproof} (of Theorem \ref{thm:s1_universality})
\begin{enumerate}
\item Observe that $[0,T)$ is a fundamental domain for $u$.  Also note that $\phi(x) = x/T \mod 1$ is a surjective submersion.  Hence with $U(\theta)= u(\theta/T)$ for $0 \le \theta \le 1$, then $u=U \circ \phi$. 

\item Suppose that $u$ has another quasiperiodic factorization $u=U' \circ \phi'$ so that $\phi':\mathbb{R} \to C$ is a surjective submersion.  Because $u$ is not constant and $\phi'$ is a submersion, $\text{rank }\phi' = 1$.  Therefore $\phi'$ has no critical points.  Since $\phi'$ is surjective, we observe that $C$ must be a manifold of dimension 1 without boundary.  Therefore, $C$ is of two possible types:
\begin{enumerate}
\item $C$ is diffeomorphic to $(a,b) \subseteq \mathbb{R}$.  Because $\phi'$ has no critical points, it must be monotonic and therefore injective.  Since it is a submersion, $\phi'$ has a continuous inverse by the inverse function theorem.  Construct a map $c:(a,b)\to S^1$ by $c(x) = \phi'^{-1}(x)/T \mod 1$.  This map is the required one, since
\begin{eqnarray*}
c(\phi'(y)) &=& \phi'^{-1}\left(\phi'(y)\right)/T \mod 1 \\
&=& y/T \mod 1 = \phi(y).\\
\end{eqnarray*}
\item $C$ is diffeomorphic to $S^1$.  In this case, we can simply write $\phi':\mathbb{R} \to S^1$.  By hypothesis, $\phi'$ is injective on every interval of the form $[a,a+T)$.  Since $u=U'\circ \phi'$ is periodic with period $T$, $\phi'$ must also be periodic with period a multiple of $T$.  Thus we simply define $c:S^1\to S^1$ by $c(x) = \phi'^{-1}(x)/T \mod 1$.
\end{enumerate}
\end{enumerate}
\end{IEEEproof}

\section{Existence of universal quasiperiodic factorizations}
\label{sec:general_universality}

\begin{lemma}
If $u=U\circ \phi = U' \circ \phi'$ are two universal quasiperiodic factorizations with phase spaces $C$ and $C'$, respectively, then $C$ and $C'$ are diffeomorphic.
\end{lemma}
\begin{IEEEproof}
Observe that because the factorizations are both universal, there are smooth maps $c:C\to C'$ and $c':C' \to C$ making the diagrams below commute.
\begin{equation*}
\xymatrix{
M\ar[r]^\phi \ar[d]_{\phi'}& C\ar[d]^U &   M\ar[r]^\phi\ar[d]_{\phi'} & C\ar[d]^U \ar[dl]^c \\
C'\ar[r]_{U'}\ar[ru]^{c'}&N           &   C'\ar[r]_{U'} & N \\
}
\end{equation*}
because both $\phi$ and $\phi'$ are surjective, this implies that both $c$ and $c'$ are surjective.  

Suppose $x,y \in C$ such that $c(x) = c(y) \in C'$.  Since $\phi$ is surjective, there are $x',y' \in M$ so that $\phi(x')=x$ and $\phi(y')=y$.  Since the diagrams above commute, $c' \circ \phi' = \phi$, which means that $c'(\phi'(x'))=x$ and $c'(\phi'(y'))=y$.  However, $\phi'(x')=c(x)=c(y)=\phi'(y')$, so we must conclude that $x=y$.  Thus both $c$ and $c'$ are bijective and are each others' inverses.  Since they are both smooth, they are diffeomorphisms.
\end{IEEEproof}

\begin{lemma}
\label{lem:pairwise_factor}
Suppose that $u=U\circ \phi = U' \circ \phi'$ are both quasiperiodic factorizations of the same function $u:M\to N$ with phase spaces $C$, $C'$ respectively.  Then there is a quasiperiodic factorization $u=U'' \circ \phi''$ and maps $p,p'$ such that the following diagram commutes:
\begin{equation*}
\xymatrix{
&C\ar[d]^p\ar[dr]^U&\\
M\ar[ru]^\phi \ar[r]^{\phi''} \ar[rd]_{\phi'}&C''\ar[r]^{U''}&N\\
&C'\ar[u]^{p'}\ar[ru]_{U'}&\\
}
\end{equation*}
\end{lemma}

\begin{IEEEproof}
Consider $C''=(C \sqcup C') / \sim$ where $x\in C \sim y\in C'$ if there is a $z \in M$ such that $\phi(z) = x$ and $\phi'(z)=y$.  Let $p:C\to C''$ and $p':C' \to C''$ be the canonical projections.  Then the function $\phi'':M\to C''$ is given by $\phi'' = p \circ \phi = p' \circ \phi'$.  It is well-defined by construction: if $z\in M$, then $\phi(z) \sim \phi'(z)$ which have the same images in $C''$ through $p$ and $p'$.

$\phi''$ is surjective because both $\phi$ and $\phi'$ are, since if $[x]\in C''$ for some $x\in C$, then surjectivity of $\phi$ implies that there is a $z\in M$ with $\phi(z)=x$.  

What remains to be shown is that $C''$ is a manifold and that $\phi''$ is a submersion.  To this end, let $[x]\in C''$.  Without loss of generality, assume that $x,y\in [x]$ for some $x\in C$ and $y\in C'$.  Let $z\in \phi^{-1}(x)$.  By the constant rank theorem (\cite[Thm 7.13]{Lee_2003}), the tangent plane at $z$ splits
\begin{equation*}
T_z M = (d_z\phi)^{-1}(T_x C)\oplus K = (d_z \phi')^{-1}(T_yC')\oplus K',
\end{equation*}
where $K,K'$ lie in the kernel of the derivatives of $\phi$ and $\phi'$.  Therefore, only vectors in $(d_z\phi)^{-1}(T_x C)\cap (d_z\phi')^{-1}(T_y C')$ correspond to the tangent plane to $C''$ at $[x]$ -- hence $C''$ is a manifold.  Any other vector in $T_z M$ corresponds to a direction that does not displace either $x\in C$ or its corresponding $y \in C'$.  Such points are all attached to $[x] \in C''$.  This also explains that $p$ is a surjective map from $T_x C$ to $(d_z\phi)^{-1}(T_x C)\cap (d_z\phi')^{-1}(T_y C')$.  Hence $\phi''$ is a submersion.
\end{IEEEproof}

\begin{IEEEproof}(of Theorem \ref{thm:universality_exists})
Suppose that $u:M\to N$ is a quasiperiodic function.  Let $I$ be the set of all quasiperiodic factorizations of $u$ so that for each $i \in I$, $u = U_i \circ \phi_i$ with $\phi_i: M \to C_i$ and $U_i: C_i \to N$.  This set has a partial order $\le$ in which $i \le j$ if there is a map $c_{ij}: C_i \to C_j$ such that the diagram below commutes
\begin{equation*}
\xymatrix{
M \ar[d]_{\phi_i} \ar[r]^{\phi_j} & C_j \ar[d]^{U_j} \\
C_i \ar[r]_{U_i} \ar[ur]^{c_{ij}} & N\\
}
\end{equation*}
Notice that $\le$ is transitive because if $i \le j$ and $j \le k$, then $c_{ik} = c_{jk} \circ c_{ij}$, so $i \le k$.  In this way, we construct the phase space for the universal factorization by 
\begin{equation*}
C_* = \lim_\to C_i = \bigsqcup_{i\in I} C_i / \sim
\end{equation*}
where $x_i \in C_i \sim x_j \in C_j$ if there is some $k\in I$ such that $c_{ik}(x_i) = c_{jk}(x_j)$.  

The function $\phi_*:M \to C_*$ can be constructed by taking the composition of any $\phi_i$ with the projection $C_i \to C_*$, and the function $U_*:C_* \to N$ is given by the constraint that $u = U_* \circ \phi_*$.  Notice that $U_*$ is well-defined by construction, since points in $C_i$ and $C_j$ are identified only when they have the same image in $N$.

Finally, we note that the proof of Lemma \ref{lem:pairwise_factor} easily lifts from a pair of quasiperiodic factorizations to an arbitrary collection of factorizations.  Thus $C_*$ is a manifold and $\phi_*$ is a surjective submersion.
\end{IEEEproof}

\section{Practical factorization algorithms}

For $S^1$-quasiperiodic functions that arise from the composition of a diffeomorphism and a periodic function on $\mathbb{R}$ (\emph{not} the case of Example \ref{eg:not_even_periodic}), DeSilva and his collaborators present an effective factorization algorithm using persistent cohomology \cite{DeSilva_2011}.  Their algorithm is based on a delay-coordinate embedding approach (see \cite{Sauer_1991, Takens_1981} for instance).  Suppose that $u:\mathbb{R}^n \to N$ is a quasiperiodic function for which we would like to obtain a quasiperiodic factorization.  Select a collection of $k$ vectors, $\{v_1, \dotsc, v_k\}$ and define the function $\phi:\mathbb{R}^n \to N^{k+1}$ by
\begin{equation*}
\phi(x)=(u(x),u(x+v_1),\dotsc,u(x+v_k)).
\end{equation*}
and the function $U:N^{k+1}\to N$ by $U(x_0,\dotsc,x_k)=x_0$.  Clearly, this is a factorization $u=U\circ \phi$.  According to a generalization of Takens' theorem \cite[Thm 2.7 and Rem 2.9]{Sauer_1991}, a generic collection of $\{v_i\}$ yields a $\phi$ which is an immersion if $k \ge 2n$.  If we restrict consideration to $\phi(\mathbb{R}^n)\subseteq N^{k+1}$, then usually $\phi$ is a surjective submersion.  Therefore, this approach yields a quasiperiodic factorization of $u$.

Although this approach yields quasiperiodic factorizations, they are often quite far from being \emph{universal}.  For instance, consider the the $S^1$-quasiperiodic function $u(t)=\sin\left(1/t\right)$ for $t>0$.  With $k\ge 3$, the delay-embedding approach will yield a quasiperiodic factorization, but this factorization will have $\mathbb{R}$ for a phase space, not $S^1$ as it should be if it were universal.  

Other approaches for identifying phase spaces have been suggested as well, such as the use of successive derivatives $\left(u(x),\frac{du}{dx}(x),\frac{d^2u}{dx^2}(x),\dotsc\right)$ in \cite{Packard_1980}, but these too tend to recover $\mathbb{R}$ as a phase space.

\section{Conclusion}
This paper proposes a new definition for quasiperiodic factorizations of functions and proves universal factorizations exist.  It characterizes these universal factorizations in the case of periodic functions.  Finally, it offers an algorithmic approach for producing quasiperiodic factorizations, which may fail to be universal.  It remains to discover factorization algorithms that ensure universality.



\bibliographystyle{IEEEtran}
\bibliography{sampta2015_bib}
\end{document}